\documentclass[12pt,a4paper]{amsbook}
\usepackage[latin1]{inputenc}
\usepackage{amsmath,amssymb}
\usepackage[francais]{babel}

\newtheorem{theo}{Theorem\ }
[section]
\newtheorem{lem}[theo]{Lemma\ }
\newtheorem{coro}[theo]{Corollary\ }
\newtheorem{defi}[theo]{Definition\ }
\newtheorem{prop}[theo]{Proposition\ }

\def \N{\mathbb{N}}
\def \H{\mathbb{H}}
\def \Z{\mathbb{Z}}
\def \R{\mathbb{R}}

 \def \g{\gamma}

\def \cH{\mathcal H}
\def \cC{\mathcal C} 
\def\cE{\mathcal E}

\catcode`\=\active\def{\"u}\catcode`\=\active\def{\'e}
\catcode`\=\active\def{\^a}\catcode`\=\active\def{\`a}\catcode`\=\active\def{\c c}
\catcode`\=\active\def{\^e}\catcode`\=\active\def{\`e}\catcode`\=
\active\def{\"\i }\catcode`\=\active\def{\^\i }\catcode`\=\active\def{\^o}\catcode
`\=\active\def{\^u}\catcode`\=\active\def{\`u}\catcode`\¶=\active\def¶{\delta }
\catcode`\¹=\active\def¹{\pi }\def \g{\gamma}
\def \G{\Gamma}
\def\ds{\displaystyle}

\begin{document}     

    \centerline{\Large On the  growth of nonuniform lattices  in  pinched }
 \centerline{\Large negatively curved manifolds}
\bigskip \bigskip \centerline{}

\centerline{Fran\c{c}oise Dal'bo \, \footnote{Fran\c{c}oise Dal'bo, IRMAR, Universit de Rennes-I\@,
 Campus de Beaulieu, 35042 Rennes Cedex. mail : francoise.dalbo@univ-rennes1.fr},}
 \centerline{Marc Peign\'e \& Jean-Claude Picaud\, \footnote{Marc Peign\'e \& 
Jean-Claude Picaud LMPT, UMR 6083, Facult\'e des Sciences etTechniques, Parc de Grandmont, 37200 Tours. mail : peigne@univ-tours.fr,
jean-claude.picaud@univ-tours.fr},}\centerline{Andrea Sambusetti\, \footnote{Andrea Sambusetti       Istituto di Matematica 
G. Castelnuovo Universita' "La Sapienza" di Roma      P.le Aldo Moro 5 - 00185 Roma. mail : sambuset@mat.uniroma1.it } }
\bigskip\bigskip\bigskip\bigskip

\section{Introduction}

We study  the relation between  the exponential growth rate of volume in 
a pinched  negatively curved manifold  and the critical exponent of    its  lattices.
These objects have a long and interesting story and are closely related to the geometry and the dynamical properties of the geodesic flow of the 
manifold (see e.g. \cite{BCG}, \cite{DOP},\cite{R2} and references  therein). \\

 Throughout   this paper, $X$ will denote a  complete and simply connected Riemannian
manifold of dimension $N\geq 2$   and we will assume that $X$ has {\it pinched negative curvature}, that is its sectional curvature  $K_X$ is bounded between two negative constants  $-b^2\leq -a ^2<0$. 
A  {\it Kleinian group} of $X$ is a torsion free and discrete subgroup $\Gamma$ of $Is(X)$; then, $\Gamma $ operates freely and properly discontinuously  on  $X$ and the quotient manifold  $M:=  X/\Gamma$ has a  fundamental group which can be identified with  $\Gamma$.  The group $\Gamma$  is called a {\it lattice } when  the volume of $M$ is finite ;   the lattice is said to be {\it uniform} if $M$ is compact.

Recall that the exponential growth rate of  $X$, also known as the {\it volume entropy} of $X$, is defined as
$$\omega(X)= \limsup_{R\to+ \infty}{1\over R}\ln v_X({\bf x},R) $$
where $v_X({\bf x},R)$ is the volume of 
the  open ball $B_X({\bf x}, R)$ of $X$, centered at the point ${\bf x}$ and with radius $R$. By the triangular inequality,   this quantity does not depend on the base point ${\bf x}$ 
; furthermore,   under our pinching assumption, Bishop-Gunther's comparison  theorem 
(see \cite{GHL}) implies
\begin{equation}\label{bishopgunther}
 (N-1)a\leq \omega(X)\leq (N-1)b.
 \end{equation}

The invariant $\omega(X)$ has been intensively studied when $Is(X)$ 
admits a {\it uniform} lattice $\G$. It turns out  that, in this case, $\omega(X)$ is a true limit and equals   the topological 
entropy of the geodesic flow of the compact manifold $M$ (see \cite{M}). Furthermore, with a suitable
 normalization  on the volume of $M$, it is a complete invariant of locally symmetric   metrics on $M$ (see \cite{BCG}).
 
 The second object of our interest in this paper is the {\it Poincar series}  $P_\Gamma(s, {\bf x})$ of a  Kleinian group  $\Gamma$,  defined  by
 $$
 P_\Gamma(s,  {\bf x})=\sum_{\gamma \in \Gamma} e^{-sd( {\bf x}, \gamma  {\bf x})},
 $$
 for  ${\bf x} \in X$ and $s \in \R$.
Its abscissa of convergence, called the {\it critical exponent of $\Gamma$},  is equal to 
$$\delta(\Gamma)=\limsup_{R\to \infty}{1\over R}\ln v_\Gamma( {\bf x},R),$$
where $ v_\Gamma( {\bf x},R)$ is the cardinality of the "ball"  $ B_\Gamma( {\bf x},R)
:=\{\gamma \in \Gamma/ d({\bf x}, \gamma{\bf x})\leq R\}$  ; again,  by the triangular inequality, 
$\delta(\G)$   does not depend on  ${\bf x}$.

A  way to understand the dynamic significance of the volume entropy $\omega(X)$ and its relation with $\delta(\G)$ is
to consider the Laplace transform of the $\Gamma$-invariant volume form $dv_X$ on $X$, namely
$$I_X(s)=\int_0^{+\infty}e^{-sr}v_X({\bf x},r){\rm d}r.$$ The abscissa of convergence of $I_X(s)$ coincides with $\omega(X)$.

 By a Fubini type argument, we also have
$\displaystyle{I_X(s)={1\over s}\int_Xe^{-sd({\bf x},y)}dv_X(y) }$. If $D$  is a 
Borel fundamental
domain   for the action of $\Gamma$ on $X$, we get, by invariance of $dv_X$ :
$$s I_X(s)=\sum_{\g\in \G}\int_{\gamma D}e^{-sd({\bf x},y)} dv_X(y)=
\sum_{\g\in \G}\int_De^{-sd(\gamma^{-1}{\bf x},y)} dv_X(y)$$
which, in turns, yields :
\begin{equation}
 P_{\G}(s,{\bf x})\int_De^{-sd({\bf x},y)} dv_X(y)\leq 
s I_X(s)\leq P_{\G}(s,{\bf x})\int_De^{sd({\bf x},y)}dv_X(y)\label{ineqlaplace}\end{equation}
From the left-hand side  of (\ref{ineqlaplace}) it immediately follows that we always have
\begin{equation}
\delta(\G) \leq \omega(X).
\label{delta-leq-omega}
\end{equation}
Moreover, from the right-hand side of  (\ref{ineqlaplace}), we have $\delta(\G)=\omega(X)$ when $\Gamma$ is a uniform lattice.

In this paper we shall investigate the case where \emph{$X$ admits a non-uniform lattice} $\Gamma$. 
Let us  emphasize  that, under this assumption, if $X$ also admits a uniform lattice $\Gamma_0$ then $X$ is a symmetric space  of non compact type (and  rank $1$).
Actually, as the curvature does not vanish, the manifold $X$ is not a Riemannian product; then (by \cite{Eb}, Corollary 9.2.2), either $X$ is symmetric or the isometry group of $X$ is discrete. But, in this last case,  $\Gamma_0$ would have finite  index in $Is(X)$ (see \cite{Eb} 1.9.34) and, if $\varphi$ is a parabolic isometry of $X$, then $\varphi^n$ would belong to $\Gamma_0$ for  some  $n\geq 1$, which contradicts the fact that a uniform lattice contains only axial elements.

Somewhat surprisingly, the equality $\delta(\Gamma)=\omega(X)$ may fail for a non uniform lattice $\Gamma$ ; actually,  in the last section of this paper, we shall prove

  \begin{theo}\label{theosurface}
There exists a complete and simply connected  Riemannian surface $X$ with pinched  negative curvature which admits a  non uniform  lattice  $\G$ such that 
$$
 \delta(\G)<\omega( X). $$ \end{theo}

Our construction extends to any dimension.  To explain it, recall that  
to each cuspidal end   of the quotient manifold  $X/\Gamma$ corresponds a maximal parabolic subgroup ${\mathcal P}\subset\G$, which has a  {\it lower critical exponent}   :
$$\delta^{-}  ({\mathcal P}) =\liminf_{R\to \infty}{1\over R}\ln v_{\mathcal P}({\bf x}, R).$$
In strictly negative curvature, this  exponent is nonzero, despite   the fact that ${\mathcal P}$ is virtually nilpotent  (see \cite{Bow2}).   The key point  is that, in the variable curvature setting, $\delta^{-}  ({\mathcal P})$  may be distinct from $\delta ({\mathcal P})$,  as was  suggested a long time ago to the second author by B. Bowditch ; in contrast, it is well known that  the critical exponent of any non elementary Kleinian group always  is a true limit \cite{R1}. 
We shall show  in Section 5  that the inequality $ \omega(X)> \delta ({\G})$  may appear as soon as $\delta^{-}  ({\mathcal P})< {\delta  ({\mathcal P})/ 2}$.

On the other hand, our example induces us to  introduce a notion of  pinching for  non uniform lattices which ensures  that $\omega(X)=\delta(\Gamma)$. Namely, we say
 that $\Gamma$  is {\it  parabolically $1/2$-pinched }
if for any maximal parabolic subgroup ${\mathcal P}\subset \G$, we have
\begin{equation}\label{parabolicgap}
 {\delta ({\mathcal P}) \over \delta^- ({\mathcal P}) }\leq 2
\end{equation}
We will prove

\begin{theo}\label{parabolically1/2-pinched} Let $X$ be a complete, simply connected Riemannian manifold with pinched negative curvature. Then
 for any lattice $\Gamma \subset Is(X)$ which is parabolically $1/2$-pinched,  we have $\delta(\G)=\omega(X)$.
 \end{theo}
 
Moreover, we notice that,   under the assumptions of this theorem,  the invariant $\omega(X)$ is a true limit  ; this follows from Corollary \ref{maxen}, combined with the fact that $\delta(\G)$  is  a limit.
 
We shall see that Theorem  \ref{parabolically1/2-pinched}  covers the case  of lattices in any $1/4$-pinched negatively curved manifold (i.e. ${b^2\over a^2}\leq 4$). As far as we know, even in the classical case of Riemannian  negatively curved symmetric spaces of rank one (which are $1/4$-pinched, cp. \cite{He}), there does not exist an elementary proof of this result. 
Nevertheless, for those spaces,  the equality $\omega(X)=\delta(\Gamma)$ can be
easily  deduced from  a general and deep result of A. Eskin and C. McMullen in \cite{EM} on  lattices of affine symmetric spaces, obtained by algebraic methods. In contrast, the context of
variable negative curvature forces us to use only elementary geometric arguments.

The equality $\omega(X)=\delta(\Gamma)$ actually holds under a milder geometric assumption than 
$1/4$-pinched curvature. Namely, we will say that a  manifold $M=X/\Gamma$  has {\it asymptotically $1/4$-pinched curvature} when,  for any $\epsilon >0$, there exists a compact set $C_{\epsilon} \subset M$, such that the metric is $(  {1\over 4  +\epsilon})$-pinched on 
$M\setminus C_\epsilon$. 
A direct consequence of Theorem \ref{parabolically1/2-pinched} is 

\begin{coro}\label{asymp1/4-pince} 
 Let $X$ be a complete, simply connected Riemannian manifold with pinched negative curvature and let $\Gamma$ be 
a  lattice of $X$.  If $M:=X/\Gamma$  has   asymptotically $1/4$-pinched curvature, then
 $\delta(\G)=\omega(X)$.\end{coro}
We remark that the pinching constant ${1\over 4}$   is   optimal because, for every $\epsilon>0$,  the example we construct in Theorem \ref{theosurface} can be chosen so that the curvature is 
${1\over 4+\epsilon}$-pinched.

The paper is organized as follows. Section 2 deals with elementary  geometrical estimates  inside  horoballs.
 In Section 3, we   relate the volume growth  of balls inside a  horoball  ${\mathcal H}$  with the critical exponent  of ample   parabolic subgroups preserving ${\mathcal H}$. In section 4, we first  give an elementary proof of the equality $\omega(X)=\delta(\Gamma)$
for  $\frac14$-pinched manifolds ; this is of interest since the main idea about the behavior of a ball intersecting  a  horoball  appears clearly in the proof. The proofs of  Theorem  \ref{parabolically1/2-pinched} and Corollary \ref{asymp1/4-pince}  will follow.  Section 5  is devoted to the construction of the example of Theorem \ref{theosurface} ; this relies on pretty technical results about convex functions, postponed to the Appendix.

\ 

 \noindent {\bf Acknowledgements.} {\it We thank P. Eberlein
 who explained us the structure of $Is(X)$ and G. Courtois for various helpfull discussions about volume entropy.}

We fix here once and for all   some notation  about asymptotic behavior of functions :  
 
 \noindent  
{\bf Notations.}
We shall write
{\it
$f\stackrel{c}{\preceq}g $ (or simply $f \preceq g $) when $f(R) \leq c g( R)$ for some constant $c>0$ and $R$ large enough.
The notation
$f \stackrel{c}{\asymp}g $ (or simply $f \asymp g $) means  $f \stackrel{c}{\preceq}g \stackrel{c}{\preceq}f .$

\noindent Analogously, we whall write $f \stackrel{c}{\sim}g $ (or simply $f \sim g $) when $|f(R) -  g( R)|\leq c$ for some constant $c>0$ and $R$ large enough.

\noindent The upper and lower exponential growth rates of  a function $f $ are  denoted  by 
$\omega^+(f ) \  ( $or simpler $\omega (f ) ) $ and $\omega^{-}(f )$ respectively ; namely we have 
$$\omega^{-}(f ):={\liminf _{R\to+\infty}}  \;  {\ln f(R) \over R} 
\quad and \quad 
\omega^+(f)=\omega(f):={\limsup _{R\to+\infty}}  \;  {\ln f(R) \over R} .$$

\noindent Finally, if $f $ and $g $ are two real functions, we denote by 
$f\ast g$ the discrete convolution of $f$ with $g$, defined by
$ \displaystyle{f\ast g(R)=\sum_{n=0}^{[R]}f(n)g(R-n)} 
 $ for any $R\geq 0$.}

 \section{Radial flow and geometry of horoballs}

As the curvature is bounded from above by $-a^2<0$, we have the following classical inequality    :

\begin{lem}\label{triangles} Let $T$ be a geodesic triangle with different
 vertices ${\bf x},{\bf y},{\bf z}\in X$ and angle at ${\bf y}$ greater than $\alpha >0$. Then  there is a constant $D=D(\alpha, a)$   such that
 $$d({ \bf x},{ \bf z})\geq d({\bf x},{ \bf y})+d({\bf y},{ \bf z})-D.$$
 \end{lem}\noindent{\bf Proof.} See   \cite{CI}.\\
\hspace*{1cm}\hfill$\Box$\\

Let  $X(\infty)$  be the boundary at infinity of $X$. Fix a point $\xi$  in  $X(\infty)$ and consider  its associated  {\it radial semi-flow}, $(\psi_{\xi, t})_{t\geq 0} $ defined as follows : for any ${\bf x} \in X$,
the point  $\psi_{\xi, t}({\bf x})$  lies on the geodesic ray $[{\bf x}, \xi)$  at distance $t$ from ${\bf x}$. For any horosphere $\partial \mathcal H$ centered at $\xi$, we set ${\partial \mathcal H}(t)= \psi_{\xi,  t}({\partial \mathcal H})$, and we let $d_t$ be  the distance  induced by  $d$ on the horosphere $\partial {\mathcal H}(t)$. For any points $x, y \in {\partial \mathcal H}(t)$, we have (see \cite{HIH}) 
\begin{equation}\label{distance-horo}
 {2\over a}\sinh \Bigl({a\over 2} d(x, y)\Bigr) \leq  d_t(x, y)\leq {2\over b}\sinh \Bigl({b\over 2} d(x, y)\Bigr). 
\end{equation}
By  \cite{HIH},   the differential of the map $\psi_{\xi, t} : \partial {\cH}  \to \partial {\cH} (t)$  satisfies, for any vector $v\in T(\partial {\mathcal H})$  and any $t\geq 0$ 
\begin{equation}\label{jacobi-norm}
e^{-bt}||v||\leq ||d\psi_{\xi, t}(v)||\leq e^{-at}||v||.
\end{equation}
This readily implies the estimates 
 \begin{equation}\label{jacobi-jacobian}
e^{-b(N-1)t}\leq |Jac(\psi_{\xi,  t})|\leq e^{-a(N-1)t}.
 \end{equation}
 In particular, if $\mu_t$ is the Riemannian measure induced on $\partial {\mathcal H}(t)$ by the metric on $X$, we have, for any Borel set  $A\subset \partial {\mathcal H}$  
 \begin{equation}\label{induced-measure}
e^{-b(N-1)t}\mu_0(A)\leq \mu_t(\psi_{\xi,  t}(A))= \int _{A} |Jac(\psi_{\xi,  t})|(x)d\mu_0(x) \leq e^{-a(N-1)t}\mu_0(A).
\end{equation}
If the  points ${\bf x, y}$  belong to  the horosphere $\partial{\mathcal H}$,  we set $$\displaystyle{t_{\bf x, y} = \inf\{t\geq 0/ d_t(\psi_{\xi,  t}({\bf x}), \psi_{\xi,  t}({\bf y}) )\leq 1\}}.$$
The next  lemma, which  precises   Lemma 4 in \cite{DOP}, will   be of major importance in the following.

\begin{lem}
 \label{quasigeodesique}  There exists a constant  $c=c(a,b)>0$,  only depending  on the bounds on the curvature, such that, for any horosphere $\partial{\mathcal H}$ and any ${\bf x}, {\bf y} \in \partial{\mathcal H}$, the arc 
$\gamma_{{\bf x},{\bf  y}}$ which is the ordered union of the three geodesic segments $[{\bf x},
 \psi_{\xi, t_{\bf x, y}}({\bf x})]$, $[\ \psi_{\xi, t_{\bf x, y}}({\bf x}), \psi_{\xi, t_{\bf x, y}}({\bf y})]$ and $[ \psi_{\xi, t_{\bf x, y}}({\bf y}),{\bf y}]$ is
  a $(1,c)$-quasigeodesic. Furthermore, for any $s, t\geq 0$, we have 
$$ d(\psi_{\xi, s}({\bf x}), \psi_{\xi, t}({\bf x}))  \stackrel{c}{\sim } \varphi(s, t) $$ where $\varphi$ is the function defined on $\R_{+}\times \R_{+}$ by 
$$\varphi(s, t)=\left\{\begin{array}{l}2t_{\bf x, y}-s-t\quad {\it  when}\quad s, t \leq t_{\bf x, y}\\
\vert s-t\vert\quad  {\it  otherwise.}\end{array}\right.$$
In particular,  we have  $d({\bf x}, {\bf y})\stackrel{c}{\sim } 2t_{\bf x, y}$.

\end{lem}

\noindent {\bf Proof. } 
If $d_0({\bf x},{\bf y})\leq 1$, the arc $\gamma_{{\bf x},{\bf  y}}$ is the geodesic segment $[{\bf x},{\bf  y}]$ and the lemma is obvious  in this case.
We now  assume $d_0({\bf x},{\bf y})>1$. Let $x= \psi_{\xi, t_{\bf x, y}}({\bf x})$  and 
 $y=\psi_{\xi, t_{\bf x, y}}({\bf y})$. From the right hand side of (\ref{distance-horo}),    the distance $d(x, y)$ is bounded from below by $b':={2\over b} \sinh^{-1}{b\over 2}$.

 Let us now fix a point $\xi'$ on the boundary at infinity of   the space $\H^N_a$ of constant curvature $-a^2$, and two points $x', y'$ on the same horosphere  centered at $\xi'$,  and at  distance $b'$ each from the other     on this space  ;  comparing the triangles $x \ y \ \xi$ and $x' \ y' \ \xi'$ we deduce that 
 $\widehat { x \ y \ \xi}\leq
\widehat { x' \ y' \ \xi' }\leq {\pi \over 2}-\theta$, for some constant $\theta>0$ depending only on $a$ and $b$. Since 
$\widehat {{\bf x}\  x \ y  }\geq \pi/2$, we have $\widehat {x\ y\ {\bf x}  }\leq \pi/2$ and so
$\widehat {{\bf x}\ y\ {\bf y}  }\geq \theta$. Applying Lemma  \ref{triangles} successively to the triangles  
${\bf x} \ x\ y$ (with $\alpha \geq \pi/2$) and ${\bf x}\ y \ {\bf y}$ (with $\alpha \geq \theta$) we obtain    $d({\bf x},{\bf y})\sim d({\bf x}, x)+d(y,{\bf y})$. The second point follows from  the first one, computing the distance between $\psi_{\xi, s}({\bf x}) $ and 
$\psi_{\xi, t}({\bf y})$  along $\gamma_{\bf x, y}$.
\hspace*{1cm}\hfill$\Box$

\

Applying this lemma, we  obtain  the

\begin{prop}\label{bouledanshoro}
 There exists a constant $c=c(a, b)>0$
 such that for any point $\xi$ in  $X(\infty)$, any horoball  ${\mathcal H}$
 centered at $\xi$ and any ${\bf x} \in   {\partial \cH}$  and $R>0$ we have   $$B_X(\psi_{\xi,  R/2}({\bf x}), R/2) \subset 
 B_X({\bf x}, R)\cap {\mathcal H} \subset B_X(\psi_{\xi,  R/2}({\bf x}), R/2+c).   $$
\end{prop}

\noindent {\bf Proof.} We need only to prove the second inclusion, the first one being obvious.
For ${\bf z} \in  B_X({\bf x}, R)\cap {\mathcal H}$,    denote by ${\bf y}$ the projection of ${\bf z}$ on
 ${\partial \cH}$ and by ${\bf z}_0$ the intersection of  the horosphere centered at 
$\xi$ and containing ${\bf z}$ with the geodesic ray $[\bf x,\xi)$. 

Assume first $t_{\bf x,y}\leq \mbox{max}\{ R/2, d({\bf y},{\bf z})\}$ ; setting $s=R/2$ and $t=d({\bf y}, {\bf z})$ in the previous lemma,  
 we get $d(\psi_{\xi,  R/2}({\bf x}), {\bf z})\sim |s-t|= d(\psi_{\xi,  R/2}({\bf x}), {\bf z}_0) \leq R/2$ (the last inequality following from the fact that $d({\bf x}, {\bf z}_0)\leq d({\bf x}, {\bf z})\leq R$).
 
 Assume now $t_{\bf x,y}\geq \mbox{max}\{ R/2, d({\bf y},{\bf z})\}$  ; applying twice  the previous lemma, we get in this case
$$\left\{\begin{array}{lll}d({\bf x},{\bf z})\stackrel{}{\sim } 2t_{\bf x,y}-d({\bf z},{\bf y})&
(\mbox{setting} \ s=0 \ \mbox{and}\  t=d({\bf y}, {\bf z})) &\\   
d(\psi_{\xi,  R/2}({\bf x}),{\bf z})\stackrel{}{\sim }2t_{\bf x,y}-d({\bf z},{\bf y})-R/2&
(\mbox{setting}\ s=R/2  \ \mbox{and}\   t=d({\bf y}, {\bf z})). & \end{array}\right. $$
Since ${\bf z} \in B_X({\bf x}, R)$ there, thus exists $c>0$ such that 
$d(\psi_{\xi,  R/2}({\bf x}), {\bf z})\leq R/2+c$.
\hspace*{1cm}\hfill$\Box$\\
 
In the next section, we will consider discrete parabolic subgroups  of $Is(X)$;  any such group fixes one point   $\xi \in X(\infty)$  and  preserves any horoball ${\mathcal H}$ centered at $\xi$. 
We shall investigate the relation between the critical exponent of ${\mathcal P}$ and the volume growth of $X$. Here we shall limit ourselves to remark : 

 \begin{coro}
 If $X$ is   homogeneous, then for any discrete parabolic  subgroup ${\mathcal P}$ of $Is(X)$, we have 
 $$
   \delta({\mathcal P})\leq \omega(X)/2. 
 $$ \end{coro}

 This fact is well known when $X$ is a rank one symmetric space ; Proposition   \ref{bouledanshoro} allows to understand the geometrical reason of this inequality. Actually, let ${\mathcal H}$ be an horoball preserved by ${\mathcal P}$ and let ${\bf x} \in \partial {\mathcal H}$. As
  ${\mathcal P}$ is discrete, we have $
 d:={1\over  2}\inf_{p\in {\mathcal P}}d({\bf x}, p{\bf x}) >0, 
 $
 then
 $$
 \bigsqcup_{p/d({\bf x}, p{\bf x})\leq R}B_X(p{\bf x}, d)\times [0, 1] \ \subset \ B_X({\bf x}, R+d+1)\cap {\mathcal H}.
 $$
By Proposition \ref{bouledanshoro}, we deduce
 $\displaystyle{v_{\mathcal P}({\bf x}, R) \preceq \sup_{{\bf y}\in {\mathcal H}} v_X\Bigl({\bf y}, {R+d +1\over 2}+c\Bigr)}$.  As $X$ is homogeneous, for any $\epsilon >0$, we have $v_X({\bf y}, r) \preceq e^{(\omega(X)+\epsilon)r}$ uniformly in ${\bf y}$. The Corollary follows.
 \hspace*{1cm}\hfill$\Box$\\

\section{ Growth of ample parabolic subgroups}

Let be  ${\mathcal P}$ a parabolic subgroup of $Is (X)$ fixing $\xi \in X(\infty)$.
We shall say that ${\mathcal P}$  is {\bf ample} if it   acts cocompactly on every horoball  $\partial {\mathcal H}$ centered at $\xi$. This   holds in particular when ${\mathcal P}$ is a maximal parabolic subgroup of a non uniform lattice of $Is(X)$.

We then  fix a (relatively compact) Borel fundamental domain ${\mathcal C} \subset \partial {\mathcal H}$ for the action of  ${\mathcal P}$ on ${\partial \mathcal H}$.
For any $t\geq 0$, the 
 set
${\mathcal C}_t:= \psi_{\xi, t}({\mathcal C})$ is a fundamental domain for the action of ${\mathcal P}$ on 
$\partial {\mathcal H}(t)$ ; in the same way, the set 
$\displaystyle{{\mathcal E}:= \cup_{t\geq 0}{\mathcal C}_t}$, which is canonically homeomorphic to ${\mathcal C}\times {\R}^+$,  is a fundamental domain for the action of ${\mathcal P}$ on the horoball ${\mathcal H}$.

We now associate to any ample parabolic group ${\mathcal P}$ a function ${\mathcal A}_{\mathcal P}$ which will play a crucial role in this paper : 
\begin{defi}
{\rm The {\bf horospherical area } of $\mathcal P$ 
is the function 
${\mathcal A}_{\mathcal P}({\bf x}, t)$ defined by }
$$
\forall \; {\bf x} \in \partial{\mathcal H}, \forall \; t\geq 0 \qquad {\mathcal A}_{\mathcal P}({\bf x}, t):= \mu_t(\psi_{\xi, t}({\mathcal C})).
$$
\end{defi}
The function $t\mapsto {\mathcal A}_{\mathcal P}({\bf x}, t)$   is decreasing and does not depend on the choice of the fundamental domain  ${\mathcal C}$ ; furthermore, by inequalities 
(\ref{induced-measure}),
for any $R$ and $R_0>0$, we have
\begin{equation}\label{estime-aire}
e^{-(N-1)bR_0}\leq {{\mathcal A}_{\mathcal P}({\bf x}, R+R_0) \over {\mathcal A}_{\mathcal P}({\bf x}, R)}
\leq e^{-(N-1)aR_0}.\end{equation}
The following proposition stresses the relation between the function ${\mathcal A}_{\mathcal P}$ and 
the orbital counting function $v_{\mathcal P}({\bf x}, R)$ of $\mathcal{P}$.
\begin{prop}\label{critical-exponent-of parabolic}There exists a constant $c=c(a, b, diam({\mathcal C}))>0$ such that for any ${\bf x} \in X$ 
$$v_{\mathcal P}({\bf x}, R)\stackrel{c}{\asymp}
{1\over
{\mathcal A}_{\mathcal P}({\bf x}, {R \over 2})
}.$$
In particular, we have 
\begin{equation}  \delta({\mathcal P})= \omega\Bigl({1\over
{\mathcal A}_{\mathcal P}({\bf x}, {R\over  2})}\Bigr) \quad \mbox{\it and }\quad 
\delta^-({\mathcal P})= \omega^-\Bigl({1\over
{\mathcal A}_{\mathcal P}({\bf x}, {R\over  2})}\Bigr).
\end{equation}
\end{prop}

\noindent {\bf Proof.}
We recall that   $d_{t}$  denotes the horospherical distance on the   horosphere
$\partial {\mathcal H}(t)$. We let $c$ be the constant of Lemma \ref{quasigeodesique} such that
$d({\bf x}, {\bf y})\stackrel{c}{\sim}2t_{{\bf x}, {\bf y}} $ for ${\bf x},   {\bf y}$  on 
$\partial {\mathcal H}$. If $d({\bf x}, {\bf y})=R$, as $t_{{\bf x}, {\bf y}} \stackrel{c/2}{\sim}{R \over 2}$, we deduce
$$
d_{R+c\over 2}\Bigl({\bf x}({R+c\over 2}), {\bf y}({R+c \over 2})\Bigr)\leq 1 \quad \mbox{\rm and}\quad 
d_{R-c\over 2}\Bigl({\bf x}({R-c \over 2}), {\bf y}({R-c \over 2})\Bigr)\geq 1.
$$
This implies that 
$\displaystyle{
\psi_{ {R+c\over 2}}(B_X({\bf x}, R)\cap \partial {\mathcal H} )\subset B_1
}$ and $\displaystyle{
\psi_{  {R-c\over 2}}(B_X({\bf x}, R)\cap \partial {\mathcal H} )\subset  B_2}$
 with $$B_1:=B_{\partial{\mathcal H}({R+c\over 2})}\Bigl( {\bf x}({R+c \over 2}), 1\Bigr)\quad \mbox{\rm and} \quad 
B_2:= B_{\partial{\mathcal H}({R-c\over 2})}\Bigl( {\bf x}({R-c \over 2}), 1\Bigr).$$
 Gauss equation implies that the sectional curvature of all horospheres for the induced  metric  is in between $ a^2-b^2 $ and $ 2b(b-a)$ (see  (\cite{BK}, section $ 1.4$, example
$(iii) $).
Therefore, there exist positive constants $v^-=v^-( a, b, {\bf x})$ and $v^+=v^+( a, b, {\bf x})$  such that  $v^- \leq vol(B_i)\leq v^+ $  for the induced volume form on the horospheres and $i=1,2$.

Now, there are  at most $v_{\mathcal P}({\bf x}, R)$ distinct fundamental domains $p({\mathcal C})$ included in $B_X({\bf x}, R)\cap \partial {\mathcal H}$ and  since the radial semi-flow $(\psi_{\xi,  t})_{t\geq 0}$ is equivariant with respect to the action of ${\mathcal P}$ on the horospheres $\partial {\mathcal H}(t)$, there are also at most $v_{\mathcal P}({\bf x}, R)$ distinct fundamental domains $p({\mathcal C}({R+c \over 2}))$ included in $\displaystyle{
\psi_{  {R+c\over 2}}(B_X({\bf x}, R)\cap \partial {\mathcal H} )}$. Therefore, we have
$\displaystyle{
v_{\mathcal P}({\bf x}, R)\leq {v^+\over {\mathcal A}_{\mathcal P}({\bf x}, {R+c\over 2})}
}$ and by (\ref{estime-aire}), this leads to 
$$
v_{\mathcal P}({\bf x}, R )\preceq{1\over {\mathcal A}_{\mathcal P}({\bf x}, {R \over 2})}.
$$
On the other hand, we can cover the set $B_X({\bf x}, R)\cap\partial {\mathcal H}$ with 
$v_{\mathcal P}({\bf x}, R+d)$ distinct fundamental domains $p({\mathcal C})$ ; by the equivariance of $(\psi_{\xi,  t})_t$ we deduce again that 
$\displaystyle{
\psi_{ {R-c\over 2}}(B_X({\bf x}, R)\cap {\mathcal H} )
}$ can be covered by $ v_{\mathcal P}({\bf x}, R+d))$ fundamental domains as well. Therefore, using 
(\ref{estime-aire}) again
$$v_{\mathcal P}({\bf x}, R )\geq {v^-\over {\mathcal A}_{\mathcal P}({\bf x}, {R-c-d\over 2})}
\succeq
  {1\over {\mathcal A}_{\mathcal P}({\bf x}, {R \over 2})}.
$$
\hspace*{1cm}\hfill$\Box$\\
We now estimate the volume of  a ball of radius $R$, inside the horoball ${\mathcal H}$. We have
\begin{prop}\label{ball-inter-horo}
There exists a constant $c=c(a, b, diam({\mathcal C}))>0$ such that
$$vol(B_X({\bf x}, R)\cap {\mathcal H})\stackrel{c}{\asymp}\int_0^{R}{ {\mathcal A}_{\mathcal P}({\bf x}, t)
\over  {\mathcal A}_{\mathcal P}({\bf x},{t+R\over 2})}dt.$$
\end{prop}

To get this result, we need the following  refinement of  Proposition \ref{bouledanshoro}.

\begin{lem}\label{hauteurs}
There exists a constant $\Delta=\Delta(a, b, diam({\mathcal C}))$ such that 
$$p ({\mathcal C})\times \Bigl[(2t_p-R+\Delta)^+,(R-\Delta)^+\Bigr[\subset\Bigl( p ({\mathcal E})\cap B_X({\bf x},R)\Bigr)\subset p ({\mathcal C})\times \Bigl[(2t_p-R-\Delta)^+,R\Bigr].$$
 \end{lem}

\noindent {\bf Proof.} Let  $\Delta=c +diam({\mathcal C}) ,$ where $c$ is the constant of Lemma 
\ref{quasigeodesique}.
We first prove the right hand side inclusion.
Let ${\bf z}= ({\bf z}_0, t) \in p({\mathcal C})\times \R^+$ and assume that this point belongs to $B_X({\bf x}, R)$.
Clearly $t\leq R$ as $t=B_\xi({\bf x}, {\bf z})\leq d({\bf x}, {\bf z})\leq R$.  If $t_p \leq {R+\Delta\over 2}$ there is nothing left to prove ; on the other hand, if $t_p>{R+\Delta\over 2}$, then $2t_p-t\stackrel{c}{\sim} d({\bf x}, {\bf z})<R$ hence $t \in [(2t_p-R+\Delta)^+, R]$. Let us now consider the case where 
${\bf z} \in  p({\mathcal C})\times  [(2t_p-R+\Delta)^+, (R-\Delta)^+]$. We may assume that $R\geq \Delta$ and $t_p\leq R-\Delta$, otherwise there is nothing to prove. If $t\geq t_p$ we have
$d({\bf x}, {\bf z})\stackrel{c}{\sim}t\leq R-\Delta$, otherwise we have 
$d({\bf x}, {\bf z})\stackrel{c}{\sim}2t_p-t\leq 2t_p-(2t_p-R+\Delta)^+$ ; therefore, in both cases ${\bf z}\in B_X({\bf x}, R)$.
\hspace*{1cm}\hfill$\Box$\\

\noindent {\bf Proof of Proposition \ref{ball-inter-horo}
}. We  simply write ${\mathcal A}(R)= {\mathcal A}_{\mathcal P}({\bf x}, R)$. Recall that
$$B_X({\bf x}, R)\cap {\mathcal H}= \bigsqcup_{p\in {\mathcal P}} B_X({\bf x}, R)\cap p({\mathcal E}).$$ 
By Lemma \ref{hauteurs}, we have
$\  B_X({\bf x}, R)\cap p({\mathcal E})\subset p({\mathcal C})\times[(2t_p-R-\Delta)^+, R].$ 
Then, we find
\begin{eqnarray*}
\sum_{p\in {\mathcal P}} vol\Bigl(
B_X({\bf x}, R)\cap p({\mathcal E})
\Bigr)&=&
\sum_{  t_p\leq R+{\Delta\over 2}}
\int_{(2t_p-R-\Delta)^+}^{R }{\mathcal A}(t)dt\\
&=&\sum_{  t_p\leq R+{\Delta\over 2}}
\int_{0}^{R}{\mathcal A}(t)1_{[(2t_p-R-\Delta)^+, +\infty[}(t)dt\\
\end{eqnarray*}
Now,  as $d({\bf x}, p{\bf x})\stackrel{c}{\sim}2t_p \leq c\leq \Delta$,  for every fixed $t \in [0, R]$  we have 
\begin{eqnarray*}
\sharp\Bigl\{ p\in {\mathcal P}/t_p\leq  R+{\Delta \over 2}\ \mbox{\rm and} \  2t_p-R-\Delta \leq t \Bigr\}&\leq&
v_{\mathcal P}\Bigl({\bf x}, {t+R+\Delta\over 2}+\Delta\Bigr)\\
&\leq& {v^+\over {\mathcal A}({t+R+3\Delta\over 2})}\\
&\preceq& {1\over {\mathcal A}({R+t \over 2})},
\end{eqnarray*}
where we have successively used Proposition \ref{critical-exponent-of parabolic} and (\ref{estime-aire}). This  yields 
$$
 vol(B_X({\bf x}, R)\cap {\mathcal H})
 \preceq  
  \int_0^{R } {{\mathcal A}(t) \over {\mathcal A}({ t+R \over 2})} dt.$$
We now prove the converse inequality. Again, by Proposition \ref{hauteurs}, we deduce
$$
B_X({\bf x}, R)\cap p({\mathcal E})\supset p({\mathcal C})\times[(2t_p-R+\Delta)^+, R-\Delta].
$$ We  only consider those $p$'s such that 
$ {R-\Delta\over 2}\leq t_p\leq {R-\Delta }$ ; summing over these $p$'s, we find
\begin{eqnarray*}
\sum_{ {R-\Delta\over 2}\leq t_p\leq R-\Delta}
vol\Bigl(
B_X({\bf x}, R)\cap p({\mathcal E})\Bigr)&=&
\sum_{ {R-\Delta\over 2}\leq t_p\leq R-\Delta}\int_{2t_p-R-\Delta}^{R-\Delta}{\mathcal A}(t)dt\\
&\geq&\sum_{ {R-\Delta\over 2}\leq t_p\leq R-\Delta}\int_{R_0}^{R-\Delta}{\mathcal A}(t)1_{[2t_p-R+\Delta, R-\Delta]}(t)dt\\
\end{eqnarray*}
for any  $R_0\geq 0$.
Now,  for every fixed $t \in [R_0, R-\Delta]$, we have 
\begin{eqnarray*}
\sharp\Bigl\{ p\in {\mathcal P}/{R-\Delta\over 2}\leq t_p\leq  R-\Delta \ \mbox{\rm and}\  2t_p-R+\Delta \leq t \Bigr\}
&\geq &
v_{\mathcal P}\Bigl({\bf x}, t+R-2\Delta\Bigr)-v_{\mathcal P}\Bigl({\bf x}, R\Bigr)\\
&\geq&  { v^-\over {\mathcal A}({t+R-2\Delta\over 2})}-{v^+\over {\mathcal A}({R\over 2})}\\
&\geq&{1 \over  {\mathcal A}({t+R\over 2})}  \Bigl( v^-{ {\mathcal A}({t+R\over 2})\over {\mathcal A}({t+R-2\Delta\over 2})}-v^+{ {\mathcal A}({t+R\over 2})\over {\mathcal A}({R\over 2})}\Bigr)\\
&\geq&{1 \over  {\mathcal A}({t+R\over 2})}  \Bigl(v^- e^{-b(N-1)\Delta}-v^+e^{-a(N-1)R_0/2}\Bigr)
\end{eqnarray*}
by  Proposition \ref{critical-exponent-of parabolic} and (\ref{estime-aire}). Therefore, if $R_0$ is large enough,    we find
$$
 vol\Bigl(B_X({\bf x}, R)\cap {\mathcal H}\Bigr)  \succeq 
  \int_{R_0}^{R-\Delta}{{\mathcal A}(t) \over {\mathcal A}({ t+R \over 2})}dt. 
$$
We can replace this last integral by 
$\displaystyle{\int_{0}^{R }{{\mathcal A}(t) \over {\mathcal A}({ t+R \over 2})}dt}$ since, $\displaystyle{\int_{R-\Delta}^{R }{{\mathcal A}(t) \over {\mathcal A}({ t+R \over 2})}dt}$ is bounded in terms of $a, b$ and $\Delta$ and for $R$ large enough
$$
\int_{R_0}^{R-\Delta}{{\mathcal A}(t) \over {\mathcal A}({ t+R \over 2})} dt\geq 
\int_{R_0}^{2R_0}{{\mathcal A}(t) \over {\mathcal A}({ t+R \over 2})} dt\asymp {1\over {\mathcal A}(R/2)}
\asymp 
\int_{0}^{R_0}{{\mathcal A}(t) \over {\mathcal A}({ t+R \over 2})} dt.
$$

\hspace*{1cm}\hfill$\Box$\\
As a direct consequence of Propositions \ref{critical-exponent-of parabolic} and \ref{ball-inter-horo},
 we obtain  
 
 \begin{coro}
 For any $\epsilon >0$ and ${\bf x} \in \partial{\mathcal H}$, we have
 
  i)  if $\delta({\mathcal P})\geq 2\delta^-({\mathcal P})$ then 
 $$
 e^{(\delta^-({\mathcal P})-\epsilon)R}\preceq vol\Bigl(B_X({\bf x}, R)\cap {\mathcal H})\Bigr)
 \preceq e^{2\Bigl(\delta({\mathcal P})-\delta^-({\mathcal P})+\epsilon\Bigr)R}
$$

ii)  if $\delta({\mathcal P})< 2\delta^-({\mathcal P})$ then 
$$
 e^{(\delta^-({\mathcal P})-\epsilon)R}\preceq vol\Bigl(B_X({\bf x}, R)\cap {\mathcal H})\Bigr)
 \preceq e^{2\Bigl(\delta({\mathcal P})+\epsilon\Bigr)R}.$$
 \end{coro}

\section{ Growth of nonuniform lattices }

We suppose now that the manifold $X$ admits a nonuniform lattice $\Gamma$.
 Let us recall some well known geometrical properties of $\Gamma$ proved  in the general context of geometrically finite groups  in (\cite{Bow1}). 
 Since  the volume of $M=X/\Gamma$ is finite, the {\it limit set} of $\Gamma$ equals $X(\infty)$ and is
 the disjoint union of its  {\it radial} subset  and of finitely many orbits
$\Gamma \xi_{1}, \ldots, \Gamma \xi_{l}$ 
of points, called {\it bounded parabolic fixed points}.  By definition,  a point $\xi_i$ corresponds to  a end of the manifold 
$M$ and  is fixed by a  parabolic subgroup of $\Gamma$.  Denote ${\mathcal P}_i$ the maximal parabolic subgroup fixing the point $\xi_i$. This group preserves  any horoball ${\mathcal H}$ centered at $\xi_i$ and acts cocompactly on the horosphere $\partial {\mathcal H}$. By Margulis' lemma (see \cite{R2}), there exist  closed horoballs ${\mathcal H}_{\xi_{1}}, \ldots, {\mathcal H}_{\xi_{l}}$ centered respectively 
at $\xi_{1}, \ldots, \xi_{l}$, such that all the horoballs $\gamma . {\mathcal H}_{\xi_{i}}$, for  $1\leq i\leq l$ and  $\gamma \in \Gamma$, are 
disjoint   or coincide. 
 We fix an origin ${\bf o} \in X$ and   a convex  Borel fundamental domain ${\mathcal D}$ in $X$ for the action of $\Gamma$,  containing the geodesic rays 
$[{\bf o}, \xi_1)..., [{\bf o}, \xi_l)$. For each $1\leq i \leq l$, we set ${\mathcal E}_{i}={\mathcal D}\cap {\mathcal H}_{\xi_i}$ and ${\mathcal C}_{i}={\mathcal D}\cap  \partial {\mathcal H}_{\xi_{i}} $.  Those both sets 
 are fundamental domains for   
  the  action of the group ${\mathcal P}_{i}$   
respectively on  $ {\mathcal H}_{\xi_{i}}$ and $\partial {\mathcal H}_{\xi_{i}}$.  Moreover, the set ${\cC}_{0}={\mathcal D}\setminus\left(\ds\cup_{i=1}^l{\mathcal E}_i\right)$, and hence each ${\mathcal C}_{i}$, 
 is relatively compact. We may assume that ${\bf o}$ belongs to the interior of  ${\mathcal C}_0$.    
 
The quotient manifold $M$ is therefore  decomposed into the 
disjoint union of  a relatively compact set $ C_{0}$ and finitely many ends  of finite volume $E_{i}={\mathcal H}_{\xi_{i}}/{\mathcal P}_i$, which are the projections on $M$  of the domains  ${\cC}_0$ and $ {\cE}_i$ respectively.

We first precise some  bounds  on the critical exponent $\delta(\G)$ in terms of bounds on the curvature of $X$.

 \begin{lem} \label{ineqexposant} We have  $(N-1)a\leq\delta(\G) \leq (N-1)b.$\end{lem}
In particular, when
 $X$ is the real hyperbolic space $\H^N_a$ of constant curvature $-a^2$,  we have  $\delta(\G) = (N-1)a$ and hence $\delta(\G)= \omega(\H^N_a)$. 

\

  \noindent{\bf Proof.} The inequality $ \delta(\G)\leq (N-1)b$ follows from (\ref{delta-leq-omega}) and ( \ref{bishopgunther}). It remains to prove the left hand side inequality of the Lemma
 
 If $\delta(\G)=\omega(X)$,  the inequality  follows from( \ref{bishopgunther}).
Assume now  $\delta(\G)<\omega(X)$ and consider $s\in]\delta(\G),\omega(X)[$. Inequality  
$(\ref{ineqlaplace})$ implies
$$\int_{\mathcal D}e^{sd({\bf o},{\bf x})}dv_X({\bf x})=+\infty$$ which, by the decomposition 
${\mathcal D}={\mathcal C}_0\ds\cup\left(\ds\bigcup_{i=1}^l{\mathcal E}_i\right)$,
 is equivalent to\begin{equation}\label{cusp-infini}{\rm max}_{i\in\{1,\cdots,l\}}\int_{{\mathcal E}_i}
e^{sd({\bf o},{\bf x})}dv_X({\bf x})=+\infty.\end{equation}
Note now that for ${\bf x}\in {\mathcal E}_i$, we have  
 $B_{\xi_i}({\bf o},{\bf x})\leq d({\bf o},{\bf x}) \leq B_{\xi_i}({\bf o},{\bf x})+  diam({\mathcal C}_i)$
 where $B_{\xi_i}(.,.)$ denotes  the   
 Busemann function  centered at $\xi_i$.
 Therefore
 the integrals 
 $\displaystyle{\int_{{\mathcal E}_i}
e^{sd({\bf o},{\bf x})}dv_X({\bf x})}$ and $\displaystyle{\int_{{\mathcal E}_i}e^{sB_{\xi_i}({\bf o},{\bf x})}dv_X({\bf x})}$ are of the same nature.

By  (\ref{induced-measure}), we have $$\int_{{\mathcal E}_i}e^{sB_{\xi_i}({\bf o},{\bf x})}dv_X({\bf x})
=
\int_{d({\bf o}, {\mathcal C}_i)}^{+\infty}e^{st}\mu_t(\psi_{\xi_i,  t}({\mathcal C}_i))dt\leq \mu_0({\mathcal C}_i)
\int_0^{+\infty}e^{t[s-(N-1)a]}{\rm d}t$$
and the last
 integral must be divergent for all $s\in]\delta(\G),\omega(X)[$, so  $\delta(\G)\geq (N-1)a$. 
\hspace*{1cm}\hfill$\Box$

\

Recall that  $v_X({\bf o}, R)$ denotes the volume of the open ball $B_X({\bf o}, R)$ and that $v_\G({\bf o}, R)$ represents the cardinality of the intersection of this ball with $\Gamma(\bf o)$. The following estimate will be used to obtain a upper bound for  $\delta(\G)$.

\begin{prop}\label{ineqvolume-faible} 
There exists a constant $\Delta=\Delta(a, b, diam({\mathcal C}_0))>0$ such that, for all   $R>0$, we have  
\begin{equation}\label{ineqvol}v_X({\bf o},R-\Delta)\preceq  v_\G({\bf o},R)+\sum_{i=1}^l\sum_{n=0}^{[R]}
v_\G({\bf o},n+1)\times vol\Bigl(B_X({\bf x}_i, R-n+\Delta)\cap    {\mathcal H}_{\xi_i}\Bigr)
\end{equation}  
where ${\bf x}_i$ denotes the intersection of the geodesic ray $[\bf o, \xi_i)$ with the horosphere  $\partial {\mathcal H}_{\xi_{i}}$.
\end{prop}

\noindent {\bf Proof.}  
Set   $d_0= diam({\mathcal C}_0)$.We have 
\begin{equation}\label{balldecomposition}B_X({\bf o},R)=\Bigl( B_X({\bf o}, R)\cap \Gamma.{\mathcal C}_{0} \Bigr) \bigcup\left(\bigcup_{1\leq i\leq l}\Bigl( B_X({\bf o}, R)\cap \Gamma. {\mathcal H}_{\xi_i}\Bigr)\right)
\end{equation}
whence
$$\ds  B_X({\bf o}, R)\cap \Gamma.{\mathcal C}_{0}  \quad  \subset \bigcup_{\gamma \in B_\G({\bf o}, R+d_0)}\gamma  ({\mathcal C}_{0}) $$
and
$$vol\Bigl(B_X({\bf o}, R)\cap \Gamma.{\mathcal C}_{0}\Bigr) 
\preceq v_\Gamma(R+d_0).$$
Now, for each  $i \in \{1, ..., l\}$ 
we define a map on  $\Gamma$ 
as follows   : for any $\gamma\in \Gamma$, let 
$x_{\gamma ,i}$ be the intersection of the ray $[{\bf o}, \gamma(\xi_i))$ with the horosphere
$\gamma (\partial {\mathcal H}_{\xi_i})$. Since  ${\mathcal C}_i$ is a fundamental domain for the action of  ${\mathcal P}_{i}$ on $\partial {\mathcal H}_{\xi_i}$ there exist a finite number of elements  $\bar \gamma$ in $\gamma {\mathcal P}_{i} $ such that 
$x_{\gamma ,i}\in {\bar \gamma} ({\mathcal C}_i)$. Choose one of those elements and  denote it by  $\bar \gamma_{i} $ . Let $\bar \Gamma_{i} $ be the set of all 
$\bar \gamma_{i} $  for  $\gamma$ in $\Gamma$. Since $d(x_{\gamma ,i}, \bar \gamma_i{\bf o}) \leq d_0$, and since the angle at $x_{\gamma ,i} $ between the geodesic segments $[x_{\gamma ,i},{\bf o}]$ and $[x_{\gamma ,i},x]$
  is greater than  $\pi/2$,  by lemma \ref{triangles} there exists a constant $d_1>0$ such that 
for every $\gamma\in \G$  and $x \in \gamma{\mathcal H}_{\xi_i}\cap B_X({\bf o},R)$, we have :
$$d({\bf o},\bar\g_i{\bf o})+d(\bar\g_i{\bf o},x)-d_1\leq d({\bf o},x).$$
We  have by (\ref{balldecomposition})
$$ B_X({\bf o}, R)\cap \Gamma. {\mathcal H}_{\xi_i} 
\subset\left(  
  \bigcup_{0\leq n\leq [R+d_0]}\quad   
 \bigcup_{ \stackrel{\overline{\gamma}\in \overline{\Gamma}_i}{   n\leq d({\bf o}, \overline{\gamma}{\bf    o})
< n+1}}B_X(\overline{\gamma}{\bf o}, R-n+d_{1})\cap    \overline{\gamma}.{\mathcal H}_{\xi_i}\right).   $$
For each $i$ denote ${\bf x}_i$ the intersection of the geodesic ray $[\bf o, \xi_i)$ with the horosphere  $\partial {\mathcal H}_{\xi_{i}} $. One has  $$vol\Bigl(B_X(\overline{\gamma}{\bf o}, R-n+d_{1})\cap    \overline{\gamma}.{\mathcal H}_{\xi_i}\Bigr)
\leq 
vol\Bigl(B_X({\bf x}_i, R-n+d_{1}+d_0)\cap    {\mathcal H}_{\xi_i}\Bigr),
$$
while 
 $$\sharp \{\bar\g\in\bar\G_i/ n\leq d({\bf o},\bar\g.{\bf o})< n+1\} \leq v_\G({\bf o},n+1),$$ so
$$
 v_X({\bf o},R-d_0)\preceq  v_\G({\bf o},R)+\sum_{i=1}^l\sum_{n=0}^{[R]}
v_\G({\bf o},n+1)\times vol\Bigl(B_X({\bf x}_i, R-n+d_{1})\cap    {\mathcal H}_{\xi_i}\Bigr).
$$The lemma follows with $\Delta\geq \max(d_0, d_1)$.
\hspace*{1cm}\hfill$\Box$

\

Proposition \ref{ineqvolume-faible}   is crucial to establish Theorem \ref{parabolically1/2-pinched}  ; we first give an elementary proof of this result, in the case where $X$ is $1/4$-pinched.

\subsection{Proof of  Theorem \ref{parabolically1/2-pinched} : the ${1\over 4}-$pinched curvature case}
We prove here that if $(X,g)$ is a complete, simply connected Riemannian manifold 
with $1 / 4$-pinched negative curvature, then
 for any lattice $\Gamma \subset Is(X)$, we have $\delta(\G)=\omega(X)$.
 
 We use the notations of Proposition \ref{ineqvolume-faible}. 
  
By (\ref{delta-leq-omega}), we need only to show that  $\omega(X)\leq \delta(\Gamma)$. By Proposition \ref{bouledanshoro},  we know that for $r>0$ the set
 $ B_X({\bf x}_i, r)\cap    {\mathcal H}_{\xi_i}$ is included in the ball of radius $r/2+c$   centered at the  point $\psi_{\xi_i,  r/2}({\bf x}_i)$. Then, (\ref{ineqvol}) leads to the following inequality 
\begin{equation}
 v_X({\bf o}, R-\Delta ) \preceq   v_\G({\bf o}, R)+\sum_{n=0}^{[R]}v_\G({\bf o}, n+1)  \times
 \!\! \!\!
\sup_{{\bf x}/B_X({\bf x}, {R-n+\Delta\over 2})\subset \Theta} 
 vol\Bigl( B_X\Bigl({\bf x}, {R-n+\Delta\over 2}\Bigl)\Bigr). 
\end{equation}
 From Bishop Gunther's theorem and the fact that $\displaystyle{b^2\leq 4 a^2}$,  we have
  $$
   vol\Bigl( B_X ({\bf x},r ) \cap \Theta\Bigr) \leq  v_X({\bf x}, r)
\preceq   e^{b(N-1)r}\preceq e^{2a(N-1)r }, $$
for any ${\bf x} \in X$ and $r>0$. 
 We conclude  that 
$\omega(X)\leq (N-1)a\leq  \delta(\Gamma)$ using  Lemma \ref{ineqexposant}.
\hspace*{1cm}\hfill$\Box$

\

 \noindent{\bf Remark  - }  The above proof   uses in a crucial way Lemma \ref{ineqexposant} and it still works if we relax the pinching assumption as follows : 
 
 {\it For any $\epsilon >0$, there exists a compact set $C_\epsilon \subset M$ such that the curvature on $M\setminus C_\epsilon$ belongs to $[- (4+\epsilon) a^2, -a^2]$.} 
 
 However, this  condition is much stronger than  the $\Bigl({1\over 4+\epsilon}\Bigr)$-pinching assumption and the proof of Corollary \ref{asymp1/4-pince} requires the more  precise estimates of the volume of balls obtained in the previous section. 

\subsection{Proof of  Theorem \ref{parabolically1/2-pinched} : the general  case}

We fix here a non uniform lattice $\Gamma \subset Is(X)$ and apply the   results of Section 3  to each maximal parabolic subgroup ${\mathcal P}_i$ of $\Gamma$.  We  first set the 
\begin{defi}
{\rm Let $M=X/\Gamma$ be a complete Riemannian manifold of finite volume with $-b^2\leq K_X\leq -a^2<0$ and with  ends $E_1, ..., E_l$.  For $1\leq i\leq l$, the {\bf cuspidal function} ${\mathcal F}_i$  associated with  $E_i$ is defined by
$$
\forall {\bf x} \in X, \forall R>0 \qquad {\mathcal F}_i({\bf x}, R)=\int_{0}^{R }{{\mathcal A}_i({\bf x}, t) \over {\mathcal A}_i\Bigl({\bf x}, { t+R \over 2}\Bigr)}dt
$$
where ${\mathcal A}_i({\bf x}, t)$ is the horospherical area function associated with $E_i$.}
\end{defi}

By (\ref{estime-aire}), the growth rates $\omega^{\pm}({\mathcal F}_i({\bf x}, .))$   depend only on the  ends $E_i$ of $M$ as for any points ${\bf x, y}\in X$  and any $R_0>0$ fixed, we have 
$\displaystyle{{\mathcal F}_i({\bf x}, R) \asymp {\mathcal F}_i({\bf y}, R)}
$.
Those functions are of major importance in order to estimate $v_X({\bf x}, R)$ ; namely,  we have the

\begin{prop}\label{estime-volume}
There exists $\Delta=\Delta(a, b, diam({\mathcal C}_0))>0$ such that 

$\displaystyle{(i) \quad v_X( \bullet,  R+\Delta)\succeq v_\Gamma ( \bullet,  R)+\sum_{i=1}^{l}{\mathcal F}_i(\bullet, R)}$

$\displaystyle{(ii) \quad v_X(  \bullet, R+\Delta)\preceq v_\Gamma (\bullet, R)+\sum_{i=1}^{l}v_\Gamma(\bullet, .)\ast {\mathcal F}_i(\bullet, .)(R)}$
\end{prop}
which leads to the 
\begin{coro}\label{maxen}
We have
$\displaystyle{
\omega^\pm(X)= \max \Bigl(  \delta(\Gamma), \omega^\pm({\mathcal F}_1), ..., \omega^\pm({\mathcal F}_l)\Bigr).
}$
\end{coro}

\noindent {\bf Proof of Proposition \ref{estime-volume}}.  {\it Part (i). }
We have 
$$B_X({\bf o}, R)\supset \bigsqcup_{ \gamma \in B_\Gamma({\bf o}, R-d_0)}
\gamma ({\mathcal C}_0) \cup \bigcup_{i=1}^l\left( B_X({\bf o}, R)\cap {\mathcal H}_i\right).$$
On the other hand 
$B_X({\bf o}, R)\cap {\mathcal H}_i\supset B_X({\bf x}_i, R-d_0)\cap {\mathcal H}_i,
$ and by Proposition 
\ref{ball-inter-horo}, we have
$$
v_X({\bf o}, R)\succeq v_\Gamma({\bf o}, R-d_0)+
\sum_{i=1}^l{\mathcal F}_i({\bf x}_i, R)
$$ with ${\mathcal F}_i({\bf x}_i, R)
\asymp {\mathcal F}_i({\bf o}, R)$ ; 
the first inequality follows.

{\it Part (ii)} follows by  plugging Proposition \ref{ball-inter-horo}  in  (\ref{ineqvol}).
\hspace*{1cm}\hfill$\Box$\\

 \noindent {\bf Proof of  Theorem \ref{parabolically1/2-pinched}}
By Corollary \ref{maxen}, it is enough to show that $\omega ({\mathcal F}_i)\leq \delta(\Gamma)$ for $1\leq i \leq l$. By Proposition \ref{critical-exponent-of parabolic}, we have, for any  $\epsilon >0$ :
 $$
 {\mathcal A}_i(t)\preceq e^{-(2\delta^-({\mathcal P})-\epsilon)t}
 \quad \mbox{\rm and} \quad  
 {\mathcal A}_i\Bigl({t+R\over 2}\Bigr)\succeq e^{-( \delta({\mathcal P}) +\epsilon)(t+R)}.$$
 So, we obtain  
 $\displaystyle{{\mathcal F}_i(t)\preceq
 e^{ ( \delta({\mathcal P}) +\epsilon)R}\int_0^R e^{(\delta({\mathcal P})-2\delta^-({\mathcal P})+2\epsilon)t} dt\preceq e^{ ( \delta({\mathcal P}) +3\epsilon)R}
 }$
 as $ \delta({\mathcal P})-2\delta^-({\mathcal P})\leq 0$, therefore 
 $\omega ({\mathcal F}_i)\leq \delta({\mathcal P})\leq \delta(\Gamma)$.\hspace*{1cm}\hfill$\Box$. \\

 \noindent {\bf Proof of  Corollary \ref{asymp1/4-pince}}
Assume that $M=X/\G$ is asymptotically $\frac14$-pinched. Then, for any fixed $\epsilon >0$ we know that outside a compact subset $C_{\epsilon}$ the curvature of $M$ is between $-\beta^2$ and  $-\alpha^2$, with $\beta^2 \leq (4+\epsilon)\alpha^2$. Therefore we have
$$ e^{-\beta (N-1) t} \preceq {\mathcal A}_i (t)   \preceq e^{-\alpha (N-1) t}$$
hence, by Proposition \ref{critical-exponent-of parabolic}, we deduce that 
$${\delta ( {\mathcal P})\over \delta^-  ( {\mathcal P})  }\leq {\beta \over \alpha}\leq   2+\epsilon$$ for every maximal parabolic subgroup of $\Gamma$.
As $\epsilon$ is arbitrary, we deduce that $M$ is parabolically $\frac14$-pinched, and we conclude by  Theorem \ref{parabolically1/2-pinched}.\hspace*{1cm}\hfill$\Box$ \\

\noindent  {\bf Remark.} We have seen that, under the assumptions of Theorem \ref{parabolically1/2-pinched}, we have $\omega ({\mathcal F}_i) \leq  \delta(\Gamma)$ for all $1\leq i \leq l$ 
; in particular, 
$\omega(X)$ is  a  limit in this case.

\section{\bf An end with the leading role}

We shall construct in this section a  pinched, negatively curved surface $S=X / \Gamma$ of  finite volume  such that
$\omega(X)>\delta(\G)$.  The surface we exhibit is homeomorphic to a 3-punctured sphere, and we shall deform a hyperbolic metric  on one end $E$  of $S$. \\
Our construction rests on two main ideas:

  i) we can deform the metric  in the end $E$  varying the sectional curvature from $\alpha^2$ to $\beta^2$ on different  bands of $E$, in order that the function ${\mathcal F}$ associated to  $E$ satisfies $\omega ({\mathcal F}) > \delta ( {\mathcal P})$.

  ii)  we set $\epsilon :=  
\omega ({\mathcal F}) - \delta ( {\mathcal P})$  and  we show that   the above deformation of the metric  can be performed  in such a way that   $\delta (\Gamma) < \delta ( {\mathcal P}) + \epsilon$ also.

 By Corollary \ref{maxen}  we conclude    that $\omega(X)>\delta(\G)$.

\

 Fix positive real  numbers $\alpha$ and $\beta$ such that $\beta >2 \alpha $. We 
can construct  sequences of disjoint intervals 
$[p_n, q_n]$, $[r_n,s_n]$ included  in  $[\Delta^{n-1}, \Delta^n]$ (for some $\Delta>1)$,
 and a
$C^2$ convex,  decreasing function ${\mathcal A}(t)$
on $[\Delta,+\infty[$ whose restrictions to $[p_n, q_n]$ and $[r_n,s_n]$ coincide  respectively with $e^{-\alpha t}$ and $e^{-\beta t}$. More precisely, we can arrange the points $p_n, q_n, r_n$ and $s_n$  in order that
$q_n\geq p_n+1$ and 
$t\in[p_n,q_n]\Leftrightarrow {t+\Delta^n\over 2}\in[r_n,s_n]$, and we can choose
${\mathcal A}$ such that 
$e^{-\beta t}\leq {\mathcal A}(t)\leq e^{-\beta t}$ and ${{\mathcal A}''(t)\over {\mathcal A}(t)}\in [\alpha^2-\eta,\beta^2+\eta]$  for all $t\in [\Delta ,+\infty[$  and  some $\eta>0$.
The existence of such intervals and of the function ${\mathcal A}$ is rather technical and we postponed the details of  proof  to the Appendix (Section 6).

 By construction,   the function $\displaystyle{{\mathcal F} (R):= \int_0^R{{\mathcal A}(t)\over {\mathcal A}({t+R \over 2})}{\rm d}t}$ satisfies : 
$$\omega({\mathcal F})\geq \limsup_{n\to+\infty}{1\over \Delta^n} \ln 
\int_{p_n}^{q_n}{{\mathcal A}(t)\over {\mathcal A}({t+\Delta^n \over 2})}{\rm d}t > \beta/2.$$
 
We can now construct the surface of Theorem \ref{theosurface}.  Start from 
 a $3$-punctured sphere $S$ with a metric $g_0$ of finite volume and constant curvature $-\alpha^2$.  Let $\Gamma =\pi_1(S)$ and let ${\mathcal P}$ be the maximal parabolic subgroup associated with the end $E$ of $S$.
 Consider the horospherical parametrization 
${ \sigma}\,:\, [0,1[\times  \R^{+}\to {\mathcal E}$ of  ${\mathcal E}$ ;
with respect to these coordinates, the hyperbolic metric writes 
$g  =e^{-2\alpha t}dx^2+dt^2$. We now perturb $g $ on  $E_n=
{ \sigma}([0,1[\times [p_{n },+\infty[)$ to obtain a new $C^2$-metric $g_n$ such that  $ g_n={\mathcal A}^2(t)dx^2+dt^2$ on $E_n$,  for  ${\mathcal A}$ defined above. We shall denote by $d$ and $d_n$ the distances on $X$ associated respectively to $g$ and $g_n$ and we let 
  $\delta_n(\Gamma )$, $\delta_n({\mathcal P})$ be  the critical exponents of $\G$ and ${\mathcal P}$ relatively to the new metric $g_n$. Notice that 
$\displaystyle{K_X=-{{\mathcal A}''\over {\mathcal A}}}$ is pinched between  $-\beta^2-\eta$ and $-\alpha^2+\eta$ ; furthermore  
${\mathcal A}(R)$ is  precisely the horospherical area (length) function of ${\mathcal P}$, with respect to $g_n$, so  $\delta_n({\mathcal P}) =\beta/2$ for all $n$, by Proposition \ref{critical-exponent-of parabolic} (while $\delta^-_n({\mathcal P})\leq \alpha/2$). Since we know that $\omega ( {\mathcal F} ) = \beta/2 + \epsilon$  for some $\epsilon >0$,  it will be enough   to show that:  

\begin{prop}\label{cgap2}
 For $n $ large enough, we have $\delta_n(\G)\in]\delta_n({\mathcal P}),\delta_n({\mathcal P})+\epsilon[$.
 \end{prop}

\noindent {\bf Proof.} Let $p$ be a generator of ${\mathcal P}$ and choose another parabolic element $q\in\G$  such that $\G$ is the free non abelian group over $p$ and $q$.
Fix $N\geq 2$ ; each element  $\g\in\G\setminus\{ id\}$ can be written in a unique way as
\begin{equation}\label{decomposition}\g=p^{l_1}q^{m_1}\cdots p^{l_k}q^{m_k},\end{equation} 
where $l_i, m_i \in\Z^*$ except for $l_1$ and $m_k$ which may be  zero. Given this decomposition, we select those $l_i$  such that $|l_i|\geq N$, say $l_{i_1},\cdots l_{i_r}$, and write 
\begin{equation}\label{decomposition2}\g=Q_1p^{l_{i_1}}Q_2\cdots p^{l_{i_r}}Q_r\end{equation}
where 
 each $Q_i$ is a subword of the expression (\ref{decomposition}), containing powers 
of $q$ and powers of $p$ not exceeding $N$  in absolute
value.
Note that decomposition (\ref{decomposition2})
is still unique. We denote by ${\mathcal Q}_N$ the subset of elements $\g\in\G$ which write simply $\g=Q_1$ in (\ref{decomposition2}). 

Now let ${\bf o} \in X$ and ${\mathcal D}$ be the Dirichlet domain for the action of $\Gamma$, centered at ${\bf o}$. Roughly speaking,     the union of the geodesic segments $$[{\bf o}, Q_1({\bf o})], 
[ Q_1({\bf o}), Q_1p^{l_{i_1}}({\bf o})],\  \cdots \ , \ [ Q_1...p^{l_{i_r}}({\bf o}), \gamma({\bf o})]$$ represents a quasigeodesic which stays close to $[{\bf o}, \gamma({\bf o})]$ and each of its subsegments corresponds to the excursion of the geodesic loop $\gamma$ alternatively  outside or  inside  the cusp $E$. We now precise this argument.

As $K_X\leq -\alpha^2+\eta$, there exits a minimal angle $\theta_0>0$ such that for all ${\bf x} \in p^{\pm 2}({\mathcal D})$
and all ${\bf y} \in q^{\pm 1}({\mathcal D})$,  we have $\widehat{{\bf x} \ {\bf o}\ {\bf y}}\geq \theta_0$. Then, 
when $Q_1\neq id$ in (\ref{decomposition2}), by a ping-pong argument we deduce that 
$\angle_{\bf o}(Q_1^{-1}{\bf o},p^{l_{i_1}}Q_2\cdots Q_r{\bf o})\geq \theta_0$,  as $l_{i_1}\geq N\geq 2$.
Therefore, by  Lemma  \ref{triangles}, there exists 
 a constant $d=d(\alpha, \theta_0)>0$  such that
$$d_n({\bf o},\g ({\bf o}))\geq d_n({\bf o},Q_1({\bf o}))+d_n({\bf o},p^{l_{i_1}}Q_2\cdots p^{l_{i_{r-1}}}Q_r({\bf o}))-d$$
Repeating this argument yields 
$$d_n({\bf o},\g ({\bf o}))\geq \sum_{i=0}^rd_n({\bf o},Q_i({\bf o}))+\sum_{j=1}^{r-1}d_n({\bf o},p^{l_{i_j}}({\bf o}))-2rd.$$
 Consequently
\begin{equation}\label{majoration}\sum_{\g\in\G}e^{-sd_n({\bf o},\g ({\bf o}))}\leq \sum_{\g\in{\mathcal Q}_N}e^{-sd_n({\bf o},\g ({\bf o}))}+
\sum_{r\geq 1}\left(e^{2sd}\sum_{|k|\geq N}e^{-sd_n({\bf o},p^k({\bf o}))} \sum_{\g\in{\mathcal Q}_N}e^{-sd_n({\bf o},\g ({\bf o}))}\right)^r\end{equation}
If $n$ is  large enough
with respect to $N$, every   element of ${\mathcal Q}_N$ correspond to a geodesic loop staying in the part of $S$ where the curvature is constant
equal to $-\alpha^2$. For that choice of $n$ anf for $s={\beta +\epsilon\over 2}$, we have 
$$ \sum_{\g\in{\mathcal Q}_N}e^{-s d_n({\bf o},\g {\bf o})}\leq  \sum_{\g\in\G}e^{-s d ({\bf o},\g {\bf o})}:=A.$$
The latter series converges because
the value of the critical exponent of any lattice in the space of constant curvature case $-\alpha^2$ is $\alpha$ and $\alpha <s$. \\
Furthermore
\begin{eqnarray*} \sum_{|k|\geq N}
e^{-sd({\bf o},p^k ({\bf o}))}& \preceq &\sum_{m\geq d({\bf o},p^N( {\bf o}))}
v_{\mathcal P}({\bf o}, m)e^{-sm}\\
&\preceq& \sum_{m\geq d({\bf o},p^N ({\bf o}))}
{e^{-sm}\over {\mathcal A}({m\over 2})}\\
& \preceq& \sum_{m\geq d({\bf o},p^N ({\bf o}))}e^{-(s-{\beta\over 2})m}=  \sum_{m\geq d({\bf o},p^N ({\bf o}))}e^{-\epsilon
m/2}\end{eqnarray*}
so that $\sum_{|k|\geq N}e^{-sd({\bf o},p^k ({\bf o}))}\to 0$ when $N\to +\infty$. 
Then, we can choose $N$ and $n $   such that
$$ \sum_{\g\in{\mathcal Q}_N}e^{-s d({\bf o},\g( {\bf o}))}\leq A<+\infty  \hspace{1cm}{\rm
 and}\hspace{1cm}
\left(e^{2s d}\sum_{|k|\geq N}e^{-s d({\bf o},p^k ({\bf o}))} A\right)<1.$$
For that choice, (\ref{majoration}) implies that the Poincaré series associated with $\G$ converges at $s$ and consequently : 
$\delta(\G)\leq s<\delta(\G)+\epsilon $.\\

\noindent {\bf Remark.} Notice that the curvature of $S$ is not asymptotically ${1\over 4}$-pinched as $\beta >2 \alpha  $ ; but, letting $\alpha \to \beta/2$ and $\eta\to 0$,  the metric can be choosen so that $K_S$ is asymptotically $({1\over 4+\epsilon})$-pinched, for any $\epsilon >0$.

\hspace*{1cm}\hfill$\Box$\\

\section{Appendix}

Let $t_0,\, t_1,\, t_2,\, t_3$ be four real numbers satisfying $t_0<t_1<t_2< t_3$. Denote by $\varphi_1$ a $C^2$ convex and decreasing
function on $[t_0,t_1]$ and $\varphi_2$ a $C^2$ convex and decreasing
function on $[t_2,t_3]$. A straightforward geometric argument on epigraphs of $\varphi_1$ and $\varphi_2$ shows that the following 
inequalities :
\begin{equation}\label{interpolation1bis}
\varphi_1'(t_1)(t_2-t_1)_{\stackrel{<}{(a)}}\varphi_2(t_2)-\varphi_1(t_1)_{\stackrel{<}{(b)}} \varphi_2'(t_2)(t_2-t_1)\end{equation}
are necessary and sufficient for the existence of a $C^2$ convex decreasing function $\psi$ on $[t_0,t_3]$ such that 
$\psi_{|[t_0,t_1]}\equiv \varphi_1$ and $\psi_{|[t_2,t_3]}\equiv \varphi_2$. \\

\noindent\begin{lem}\label{convexity}
Let $\alpha,\beta$ two positive reals such that $\alpha<\beta$. \\
{\bf (I\@)} Inequalities (\ref{interpolation1bis}) are satisfied for $ \varphi_1(t)=e^{-\alpha t}$ and $ \varphi_2(t)=e^{-\beta t}$
when $t_2-t_1>{1\over \alpha}$.\\
 {\bf (II\@)} Inequalities (\ref{interpolation1bis}) are satisfied for 
$\varphi_1(t)=e^{-\beta t}$ and $ \varphi_2(t)=e^{-\alpha t}$ when $t_2>({\beta\over \alpha}+\epsilon)t_1$ for any $\epsilon>0$.\end{lem} 

\bigskip

\noindent {\bf Proof.} 
 Case {\bf (I)} : 
$$(a)\Leftrightarrow e^{-\beta t_2+\alpha t_1}+\alpha(t_2-t_1)>1$$
and the second inequality is satisfied when $t_2-t_1>{1\over \alpha}$.  Note that this condition is optimal if we want 
such an inequality to be satisfied for arbitrary large $t_1$ because with $u=t_2-t_1$, this inequality becomes
$$e^{(\alpha-\beta) t_1-\beta u}+\alpha u>1$$
and this inequality cant be satisfied for small $u$ when $t_1$ is too large.\\
With the previous notations, 
$$(b)\Leftrightarrow e^{\beta u}e^{(\beta-\alpha) t_1}-\beta u-1>0$$
and the latter inequality is always satisfied because $e^x-x-1> 0$ for all $x>0$.\\

\noindent Case {\bf (II)} : 
 $$(a)\Leftrightarrow e^{-\alpha t_2+\beta t_1}+\beta(t_2-t_1)>1$$
and this second inequality is satisfied when $t_2-t_1>{1\over \beta}$. The same remark as in the case $(I)$.\\
With the previous notations too, 
$$(b)\Leftrightarrow e^{\alpha u}e^{(\alpha-\beta) t_1}-\alpha u-1>0$$
with $u=t_2-t_1$. If we set 
$t_2=(1+x)t_1+f(t_1)$ and substitute in the last term, a necessary condition in order to realise $(b)$ is 
$(x+1)\geq {\beta\over \alpha}$ and if we set $(x+1)= {\beta\over \alpha}$ and replace, we get
$e^{\alpha f(t_1)}-(\beta-\alpha)t_1-f(t_1)-1>0$. The conclusion follows.\\
\hspace*{1cm}\hfill$\Box$\\

\begin{lem}\label{interpolation 12}
Let $t_0<t_1<t_2<t_3$ and $\eta>0$. There exists $A=A(\eta,\alpha,\beta)>0$ and $B=B(\alpha,\beta)>0$ such that if $t_2>A.t_1$ and $t_0>B$, \\
\noindent{\bf (I)} There exists  a $C^2$ convex and decreasing function $\psi$ on $[t_0,t_3]$ satisfying :
$${\bf (C_1)}\left\{\begin{array}{llll}\forall\, t\in[t_0,t_1],&\psi(t)=e^{-\alpha t}&&\\
                                                     \forall\, t\in[t_2,t_3],&\psi(t)=e^{-\beta t}&&\\
                                                     \forall\, t\in[t_0,t_3],&\alpha^2-\eta\leq{\psi''(t)\over \psi(t)}\leq\beta^2-\eta& \mbox{and}&\psi(t)\geq e^{-\beta t}\end{array}\right.$$

\noindent{\bf (II)} There exists a $C^2$ convex and decreasing function $\psi$ on $[t_0,t_3]$ such that  we have
$${\bf (C_2)}\left\{\begin{array}{llll}\forall\, t\in[t_0,t_1],&\psi(t)=e^{-\beta t}&&\\
                                                     \forall\, t\in[t_2,t_3],&\psi(t)=e^{-\alpha t}&&\\
                                                     \forall\, t\in[t_0,t_3],&\alpha^2-\eta\leq{\psi''(t)\over \psi(t)}\leq\beta^2+\eta& \mbox{and}&\psi(t)\geq e^{-\beta t}\end{array}\right.$$

\end{lem}

\noindent {\bf Proof.} By the previous remark, if we choose  $A>{\beta\over \alpha}$ and $B>{1\over \beta-\alpha}$, 
inequalities (\ref{interpolation1bis}) are satisfied. In both cases, set
$$\psi(t)= e^{-t\varphi(t)}\;\;\; t\in[t_0,t_3]$$
where $\varphi$ is constant on $[t_0,t_1]$ and $[t_2,t_3]$ (depending in an obvious way on case I or II). 
Consider a $C^2$ function $\phi\,:\,[0,1]\rightarrow[\alpha,\beta]$ ;
 set $s={\lambda(t-t_1)}$ where $\lambda= {1\over t_2-t_1}$ and put $\varphi(t)=\phi(s)$ for $t\in[t_1,t_2]$. A straightforward calculus gives, for $s\in[0,1]$ :
$$\begin{array}{ll}\ds{\psi''(t)\over \psi(t)}&=\left( (s\phi(s))'+\lambda t_1\phi'(s)\right)^2-\lambda\left( 2\phi'(s)+(s+\lambda t_1)\phi''(s)\right)\\
                                                                   &=(k'(s))^2+\lambda t_1(2k'(s)\phi'(s)+\lambda (t_1(\phi'(s))^2-\phi''(s)))
-\lambda(2\phi'(s)+s\phi''(s))\\
&=(k'(s))^2+\theta(\lambda)\end{array}$$
where $k(s):=s\phi(s)$ and $\theta$ is a function such that $\theta(\lambda)\to 0$ when $\lambda\to 0$.\\
Set $M_i=\sup_{s\in[0,1]}|\phi^{(i)}(s)|$ for $i=1,2$ (which depend only on $(\alpha,\beta)$) and 
$C={1\over 8(\beta+1)(M_1+M_2+\beta)}$.  The  previous equalities  implie
\begin{equation}\label{controlcourbure}(k'(s))^2-{\eta\over 2}\leq {\psi''(t)\over \psi(t)}\leq (k'(s))^2+{\eta\over 2}\end{equation}
when $\lambda t_1<C.\eta$ i.e. for $t_2>(1+{1\over C.\eta})t_1:=A.t_1$. We show  in both cases that we can choose 
a $C^2$ function $\phi$ with values in $[\alpha,\beta]$ such that for all $s\in[0,1]$ :
\begin{equation}\label{controlcourbure2}\alpha-{\eta\over 4}\leq k'(s)\leq \beta+{\eta\over 4}.\end{equation}

\noindent Case {\bf (I)} : choose $\phi\,:\, [0,1]\to [\alpha,\beta]$  non decreasing satisfying $\phi(0)=\alpha$, $\phi(1)=\beta$ and
$\phi'(0)=\phi'(1)=\phi''(0)=\phi''(1)=0$. Then, the function $\varphi$ can be extend on $[t_0,t_3]$ in a $C^2$
manner and on $[0,1]$, we have $k'(s)=(s\phi(s))'=\phi(s)+s\phi'(s)\geq\alpha$ and $\phi(s)\leq \beta$ so that 
${\psi''/\psi}\geq \alpha^2-\eta$ and $\psi(t)\geq e^{-\beta t}$ are both satisfied on $[t_0,t_3]$. It implies in particular that the function
$\psi$ constructed is convex on $[t_0,t_3]$. Note that in this case, the inequality
 $\lambda.t_1<C\eta$ must be satisfied, for, in the second expression
of ${\psi''\over\psi}$, the term $(t_1(\phi'(s))^2-\phi''(s))$ is negative in the neighborhood of $s_0=\inf\{s\,;\, \phi'(s)=0\}$.\\
It is left to show that $\phi$ or equivalently $k$ can be choosen so that $k'(s)=\phi(s)+s\phi'(s)\leq \beta+{\eta\over 4}$. 
The boundary conditions for $\phi$ up to the first order translate to  $k(0)=0$, $k(1)=\beta$, $k'(0)=\alpha$ and $k'(1)=\beta$. 
For $\epsilon_1\in]0,1[$,  consider the $C^0$-piecewise affine function $\bar k$ defined on $[0,\epsilon_1]$ by $\bar k(t)=\alpha.t$, on 
$[1-\epsilon_1,1]$ by $\bar k(t)=\beta.t$ and  affine on $[\epsilon_1,1-\epsilon_1]$. 
If we choose $\epsilon_1$  small enough (depending on $\eta$ and $\alpha$), we can smooth
  $\bar k$ to obtain a $C^2$ function $k$ on $[0,1]$ 
  in such a way that the dérivative $k'$ satisfies
$$\left\{\begin{array}{ll}k'(s)=\alpha&s\in[0,\epsilon_1/2]\\
                                 k'(s)\leq \beta+\eta/(4\beta)&s\in[\epsilon_1/2,1-\epsilon_1/2]\\
                                 k'(s)=\beta&s\in[1-\epsilon_1/2,1]\end{array}\right.$$
so that $(k'(s))^2\leq \beta^2-\eta/2$.

Case {\bf (II)} : this case is similar. We  choose $\phi\,:\, [0,1]\to [\alpha,\beta]$ non increasing satisfying $\phi(0)=\beta$, $\phi(1)=\alpha$ and
$\phi'(0)=\phi'(1)=\phi''(0)=\phi''(1)=0$, or equivalently (up to the first order), we choose $k(s)=s\phi(s)$ satisfying
 $k(0)=0$, $k(1)=\int_0^1k'(s)ds=\alpha$, $k'(0)=\beta$ and $k'(1)=\alpha$. The construction is symmetric to the previous one.
In both cases, the desired inequalities : (\ref{controlcourbure2}), (\ref{controlcourbure}) and $e^{-\alpha t}\leq\psi(t)\leq e^{-\alpha t}$ are satisfied.
\hspace*{1cm}\hfill$\Box$\\

Let us now  construct the sequences of  intervals 
$[p_n,q_n]$, $[r_n,s_n]$ and the function ${\mathcal A}$
 we used in Section 4.
 Let $A>1$ and $B>0$ given by Lemma \ref{interpolation 12}. We set 
 $$\left\{\begin{array}{lll}
                                    p_n=(1-\lambda_0)\Delta^{n-1}+\lambda_0\Delta^n& and&r_n={p_n+\Delta^n\over 2}\\
                                    q_n=(1-\mu_0)\Delta^{n-1}+\mu_0\Delta^n& and&s_n={q_n+\Delta^n\over 2}\end{array}\right.$$
for $\Delta, \lambda_0$ and $\mu_0$ to be defined.

Fix $(\lambda_0,\mu_0)$ in the (nonempty) set $(]0,1[^2\cap \{(\lambda,\mu)\;;\; 1+\lambda-2A\mu>0\;\wedge\; \mu>\lambda\}).$
The polynomial function
$P(x)=2\lambda_0x^2+((2-A)-2\lambda_0-A\mu_0)x-A(1-\mu_0)$
tends  to infinity as $x\to +\infty$ ; thus, we can choose a positive real number $\Delta$ such that both inequalities
\begin{equation}\label{in1}\Delta>{2A-1+\lambda_0-2A\mu_0\over 1+\lambda_0-2A\mu_0}\end{equation}
                        \begin{equation}\label{in2}     P(q_0)>0\end{equation}
are satisfied.\\
Inequality (\ref{in1}) insures that $r_n>Aq_n$ and  inequality  (\ref{in2}) insures that $p_{n+1}>As_n$. By Lemma \ref{interpolation 12},
there exists $n_0\in\N^*$ and a $C^2$-convex and decreasing function ${\mathcal A}$ on $[\Delta^{n_0-1},+\infty[$ satisfying 
${{\mathcal A}''(t)\over {\mathcal A}(t)}\geq \alpha^2-\eta$ and ${\mathcal A}(t)\geq e^{-\beta t}$ for all $t\in [\Delta^{n_0-1},+\infty[$,
and such that for $n\geq n_0$, we have :
$$\left\{\begin{array}{ll}{\mathcal A}(t)=e^{-\alpha t}& \forall\, t\in[p_n,q_n]\\
                                      {\mathcal A}(t)=e^{-\beta t}& \forall\, t\in[r_n,s_n].\end{array}\right.$$
Note that by construction $t\in[p_n,q_n]\Leftrightarrow {t+R_n\over 2}\in[r_n,s_n]$.


\begin{thebibliography}{99}
    \bibitem{BFZ}   {\sc  Babillot M., Feres R., Zeghib A.} {\em Rigidit, groupe fondamental et dynamique},   
 Panorama et synthse. {\bf 13} SMF (2002).
 \bibitem{BGS} {\sc  Ballmann W., Gromov M., Schroeder V.} {\em Manifolds of non positive curvature},  
 Progress in Math. {\bf 61} Birkh\"auser, Boston (1985).
  \bibitem{Be} {\sc  Beardon A.F. } {\em  The exponent of convergence of Poincaré series},
Proc. London Math. Soc. (3) vol. 18 (1968), pp  461-483.
\bibitem{BCG} {\sc  Besson G., Courtois G., Gallot S.} 
{\em Entropies et rigidités des espaces localement symétriques de courbure strictement négative},
Geometric and Functional Analysis, Vol. 5, n$^\circ$ 5 (1995), pp  732-799.
\bibitem{Bow1} {\sc  Bowditch B.H.}
{\em Geometrical finiteness with variable negative  curvature}, Duke Math. J. vol.  77 (1995),  pp  229-274. 
\bibitem{Bow2} 
{\sc  Bowditch B.H.} {\em Discrete parabolic groups}, J. Diff. Geometry 8, (1993), pp  559-583.
\bibitem{BK} {\sc  Buser P.,   Karcher H.} {\em Gromov's almost flat manifolds}, Astérisque 81, S.M.F. (1981).
\bibitem{CI} {\sc  Corlette K.,   Iozzi A.} {\em Limit sets of isometry groups of exotic hyperbolic spaces}, Trans. A. M. S.vol.  3511999), n. 4, pp  1507-1530.
\bibitem{DOP} {\sc Dal'bo F.,  Otal J.P.,  Peigné M.}, 
{\em S\'eries de Poincar\'e des groupes g\'eom\'etriquement finis},
Israel Jour. Math. 118, (2000), pp 109-124.
\bibitem{E} {\sc  Eberlein P.}
{\em   Surfaces of non positive 
curvature},    Memoirs AMS {\bf 218} (1979). 
\bibitem{Eb} {\sc  Eberlein P.}
{\em   Geometry of nonpositively curved manifolds},   Chicago Lectures in Mathematics.
\bibitem{EO} 
{\sc  Eberlein P., O'Neill B.} {\em   Visibility manifolds},  Pac. J. Math.  {\bf 46} (1973) pp 45-110. 
\bibitem{EM} {\sc Eskin A.,   McMullen C.} {\em Mixing, counting and equidistribution in Lie groups}, 
 (1995) pp  181-209.
\bibitem{GHL} {\sc  Gallot S., Hulin D.,    Lafontaine J.} {\em Riemannian Geometry}, Springer-Verlag Universitext  2nd Edition.
\bibitem{He} {\sc Helgason S. } {\em Differential geometry , Lie groups and Symmetric Spaces}, Graduate Studies in mathematics  vol. 34 , (1978).
\bibitem{HIH} {\sc  Heinze  E.,   Im Hof H.C.} {\em Geometry of horospheres}, J. Diff. Geom. vol.  12, (1977) pp 481-491.
\bibitem{M} {\sc Manning A.} {\em Topological entropy for geodesic flows}, Annals of Mathematics, {\bf 110} (1979), pp 567-573.
\bibitem{OP} {\sc Otal J.P. .,  Peign\'e M. } {\em On some exotic kleinians groups}, Preprint.
\bibitem{R1} {\sc  Roblin T.} {\em Sur la fonction orbitale des groupes discrets en courbure n\'egative}, 
Ann. Inst. Fourier (Grenoble) {\bf 52} (2002) n$^{\circ}$1, pp 145-151.
\bibitem{R2} {\sc  Roblin T.} {\em Ergodicité et équidistribution en courbure négative}, Mémoires S.M.F. n.95 (2003).
\bibitem{Su} {\sc  Sullivan D.}{\em The density at infinity of a discrete group of hyperbolic motions},  Publ. Math. IHES, vol. 50, (1979), pp 171-202.
\end{thebibliography}
\end{document}